\documentclass[12pt]{article}
\usepackage{amsfonts,amssymb,epsfig}
\usepackage{amsmath}

\date{}
\begin{document}
\newtheorem{df}{Definition}
\newtheorem{thm}{Theorem}
\newtheorem{lm}{Lemma}
\newtheorem{pr}{Proposition}
\newtheorem{co}{Corollary}
\newtheorem{re}{Remark}
\newtheorem{note}{Note}
\newtheorem{claim}{Claim}
\newtheorem{problem}{Problem}

\def\R{{\mathbb R}}

\def\E{\mathbb{E}}
\def\calF{{\cal F}}
\def\N{\mathbb{N}}
\def\calN{{\cal N}}
\def\calH{{\cal H}}
\def\n{\nu}
\def\a{\alpha}
\def\d{\delta}
\def\t{\theta}
\def\e{\varepsilon}
\def\t{\theta}
\def\pf{ \noindent {\bf Proof: \  }}
\def\trace{\rm trace}
\newcommand{\qed}{\hfill\vrule height6pt
width6pt depth0pt}
\def\endpf{\qed \medskip} \def\colon{{:}\;}
\setcounter{footnote}{0}

\def\Lip{{\rm Lip}}

\renewcommand{\qed}{\hfill\vrule height6pt  width6pt depth0pt}

\title{Dimension Reduction in $L_p$, $0<p<2$
}

\author{Gideon Schechtman\thanks{Supported by the
Israel Science Foundation.\ AMS subject classification: 46B85
}}
\date{October 10, 2011}
\maketitle

\begin{abstract}

Complementing a recent observation of Newman and Rabinovich for $p=1$ we observe here that for all $0<p<2$ any $k$ points in $L_p$ embeds with distortion $(1+\e)$ into $\ell_p^n$ where $n$ is linear in $k$ (and polynomial in $\e^{-1}$).

\end{abstract}

\section{Introduction}
The very well known Johnson--Lindenstrauss Lemma \cite{jl} asserts that, for all $k$ and $\e>0$, any $k$ points in a Hilbert space embed with distortion $1+\e$ into $\ell_2^n$ for $n=O(\e^{-2}\log k)$. It is also known that nothing similar to that occurs for the $L_1$ norm: There are $k$ points in $L_1$ which if embedded in $\ell_1^n$ with distortion $D$ forces $n\ge k^{c/D^2}$ for large $D$ (\cite{bc}, and \cite{ln} for a simpler proof) and $n\ge k^{1-O(1/\log(1/(D-1)))}$ for $D$ close to 1 (\cite{acnn}). As for upper bound on $n$, until recently the best that was known was that any $n$ points in $L_p$ embed isometrically in $\ell_p^{O(k^2)}$ (\cite{b}) and with distortion $1+\e$ in $\ell_p^{O(\e^{-2}k\log k)}$ for $0<p<2$ and in $\ell_p^{O(\e^{-2}k^{p/2}\log k)}$ for $2<p<4$ (\cite{sc1}).
Recently, Newman and Rabinovich \cite{nr} observed that a recent result of Batson, Spielman and Srivastava \cite{bss} implies that one can remove the $log k$ for $p=1$ and get that $k$ points in $L_1$ $(1+\e)$-embed into $\ell_p^{O(\e^{-2}k)}$. Here we show, also using \cite{bss}, a similar result for all $0<p<2$. Our dependence on $\e$ is however worse:
Any $k$ points in $L_p$, $0<p<2$, $(1+\e)$-embed into $\ell_p^{O(\e^{-(2+\frac2p)}k)}$. The $O$ notation hides a constant depending on p.

\section{The main result}

We shall use the following two theorems. The first one is due in this form to Talagrand \cite{t} and is not known to give the best possible dependence on $\e$. (For $\rho=1/2$ there is a better result by Kahane \cite{k}.)

\begin{thm}{(\cite{t})}\label{thm:snowflake} For each $\e>0$ and $0<\rho<1$ there is a positive integer $k=k(\e,\rho)$ and a map $\varphi:\R\to\R^k$ such that
\[
(1+\e)^{-1}|x-y|^\rho\le \|\varphi(x)-\varphi(y)\|_2\le (1+\e)|x-y|^\rho
\]
for all $x,y\in\R$. Moreover, $k\le K(\rho)\e^{-1/\rho}$.
\end{thm}

The second theorem we shall use is a relatively new theorem of Batson, Spielman and Srivastava.
\begin{thm}{([BSS])}\label{thm:bss}
Suppose $0 < \e < 1$ and $A = \sum_{i=1}^mv_iv_i^T$
are given, with $v_i$ column vectors in $\R^k$. Then there are nonnegative weights $\{s_i\}_{i=1}^m$, at most $\lceil k/\e^2\rceil$
of which are nonzero, such that, putting $\tilde A=\sum_{i=1}^ms_iv_iv_i^T$,
\begin{equation}\label{eq:bss}
(1+\e)^{-2}x^TAx\le x^T\tilde Ax\le (1+\e)^{2}x^TAx
\end{equation}
for all $x\in \R^k$.
\end{thm}

We shall need the following simple corollary of this theorem which in turn is a generalization of Corollary 1 of \cite{sc2}, dealing with the case $s=1$.

\begin{co} \label{co:subspaces} Let $X_l, i=1,2,\dots,s$, be $s$ $k$-dimensional subspaces of $\ell_2^m$ and let $0<\e<1$. Then there is a set $\sigma\subset\{1,2,\dots,m\}$ of cardinality at most $n\le \e^{-2}ks$ and positive weights $\{s_i\}_{i\in\sigma}$ such that
\begin{equation}\label{eq:bssk}
(1+\e)^{-1}\|x\|_2\le (\sum_{i\in\sigma}s_ix^2(i))^{1/2}\le(1+\e)\|x\|_2
\end{equation}
for all $l=1,2,\dots,s$ and all $x=(x(1),x(2),\dots,x(m))\in X_l$.
\end{co}

\pf Let $u^l_1,u^l_2,\dots,u^l_k$ be an orthonormal basis for $X_l$, $l=1,2,\dots,s$; $u^l_j=(u^l_j(1),u^l_j(2),\dots,u^l_j(m))$, $j=1,\dots,k$.
Put $v^{lT}_i=(u^l_1(i),u^l_2(i),\dots,u^l_k(i))$, $i=1,\dots,m, \ l=1,\dots,s$. Let also $v_i$ be the concatenation of $v_i^1,v_i^2,\dots,v_i^s$ forming a column vector in $\R^{ks}$. Then $A=\sum_{i=1}^m v_iv_i^T$, is a $ks\times ks$ matrix whose $s$ $k\times k$ successive central submatrices are the $k\times k$ identity matrix. Let $s_i$ be the weights given by Theorem \ref{thm:bss} with $k$ replaces by $ks$. Let also $\sigma\subset \{1,\dots,m\}$ denote their support; The cardinality of $\sigma$ is at most $\e^{-2}ks$. Let $l=1,\dots,s$ and $x=\sum_{i=1}^k a_i u^l_i=(x(1),x(2),\dots,x(m))\in X_l$. Apply (\ref{eq:bss}) to the vector $\bar a\in \R^{ks}$ where $\bar a^T=(\bar 0,\dots,\bar 0,a^T,\bar 0,\dots,\bar 0)^T$ with $\bar 0$ denotes $0$ vector in $\R^k$ and   $a^T=(x_1,\dots,x_k)$ stand in the $(l-1)k+1$ to the $lk$ places. Then
\[
(1+\e)^{-2}\|x\|_2^2=(1+\e)^{-2}a^T\sum_{i=1}^m v^l_iv^{lT}_ia\le a^T\sum_{i=1}^m s_iv^l_iv^{lT}_ia\le (1+\e)^{2}\|x\|_2^2.
\]
Finally, notice that, for each $i=1,\dots,m$ and $l=1,\dots,s$, $a^T v^l_iv^{lT}_ia=x(i)^2$ is the square of the $i$-th coordinate of $x$. Thus,
\[
a^T\sum_{i=1}^m s_iv^l_iv^{lT}_i
a=\sum_{i=1}^ms_ix(i)^2.
\]
\endpf

The main result of this note is:
\begin{thm} For all $0<p<2$ there is a constant $K(p)$ such that for all $\e>0$ and all $z_1,z_2,\dots,z_k$ in $L_p$ there are $w_1,w_2,\dots,w_k$ in $\ell_p^n$ satisfying
\[
\|z_i-z_j\|\le\|w_i-w_j\|\le(1+\e) \|z_i-z_j\|
\]
for all $i,j$, where $n\le K(p)k/\e^{2+\frac{2}{p}}$.
\end{thm}

\pf Let $\varphi:\R\to\R^s$ with $s\le K(p)\e^{-2/p}$ be the function from Theorem \ref{thm:snowflake}:
\begin{equation}\label{eq:phi}
(1+\e)^{-1}|r-s|^{p/2}\le \|\varphi(r)-\varphi(r^\prime)\|\le (1+\e)|r-s|^{p/2}
\end{equation}
for all $r,r^\prime\in \R$. Assume as we may that $z_1,z_2,\dots,z_k\in \ell_p^m$ for some finite $m$ and consider the map $\phi:\R^m\to \R^{ms}$ given by

\[
\phi(r_1,r_2,\dots,r_m)=(\varphi(r_1),\varphi(r_2),\dots,\varphi(r_m)).
\]
Let $P_l:\R^{ms}\to\R^m$, $l=1,\dots,s$, be the restriction operator  to the coordinates $\{l,s+l,s+2l,\dots,s+(m-1)l\}$. Consider the $s$ subspaces of $\R^m$ given by
\[
X_l={\rm span}\{P_l\phi(z_1),\dots,P_l\phi(z_k)\}
\]
$l=1,\dots,s$. Apply Corollary \ref{co:subspaces} to these $s$ subspaces to get a set $\sigma\subset\{1,2,\dots,m\}$ of cardinality at most $n\le \e^{-2}ks$ and positive weights $\{s_i\}_{i\in\sigma}$ such that
\begin{equation}\label{eq:bss1}
(1+\e)^{-1}\|x\|_2\le (\sum_{i\in\sigma}s_ix^2(i))^{1/2}\le(1+\e)\|x\|_2
\end{equation}
for all $l=1,2,\dots,s$ and all $x=(x(1),x(2),\dots,x(m))\in X_l$.
Applying (\ref{eq:bss1}) to $x=P_l\phi(z_u)-P_l\phi(z_v)$ we get

\begin{align}
(1+\e)^{-2}\|P_l\phi(z_u)-P_l\phi(z_v)\|^2_2\le& \sum_{i\in\sigma}s_i(P_l\phi(z_u)-P_l\phi(z_v))^2(i)\\
\le&(1+\e)^2\|P_l\phi(z_u)-P_l\phi(z_v)\|^2_2.\nonumber
\end{align}
Adding up these $s$ inequalities, we get
\[
(1+\e)^{-2}\|\phi(z_u)-\phi(z_v)\|^2_2\le \sum_{i\in\sigma}s_i\|(\phi(z_u)-\phi(z_v))(i)\|_2^2
\le(1+\e)^2\|\phi(z_u)-\phi(z_v)\|^2_2
\]
where by $\phi(\bar r)(i)$ we mean the restriction of $\phi(\bar r)$ to the $s$ coordinates \\
$(i-1)s+1$ to $is$. Applying (\ref{eq:phi}) we now get
\[
(1+\e)^{-6}\|z_u-z_v\|^p_p\le \sum_{i\in\sigma}s_i\|(z_u-z_v)(i)\|_p^p
\le(1+\e)^6\|z_u-z_v\|^p_p
\]
for all $u$ and $v$.
\endpf


\noindent Gideon Schechtman\newline Department of
Mathematics\newline Weizmann Institute of Science\newline Rehovot,
Israel\newline E-mail: gideon.schechtman@weizmann.ac.il

\end{document}